\documentstyle[12pt]{amsart}

\newtheorem{thm}{Theorem}
\newtheorem{lem}[thm]{Lemma}

\newtheorem{pr}[thm]{Proposition}

\newcommand{\N}{{\Bbb{N}}}
\newcommand{\R}{{\Bbb{R}}}
\newcommand{\Z}{{\Bbb{Z}}}

\newcommand{\wlp}{L^{p,\infty}\msp}
\newcommand{\seq}{\el^{p,\infty}}
\newcommand{\seqr}{\el^{r,\infty}}
\newcommand{\unit}{L^{p,\infty}[0,1]}
\newcommand{\rline}{L^{p,\infty}[0,\infty)}

\newcommand{\elinf}{{\ell^\infty}}

\newcommand{\bv}{\bigvee}
\newcommand{\bc}{\bigcup}

\newcommand{\msp}{(\Omega ,\Sigma ,\mu)}
\newcommand{\set}{\{2^n : n \in \Z\}}
\newcommand{\xc}{{\{\pi\in\Pi_n : |x(\pi)| > c\}}}
\newcommand{\xgnc}{{\{(\g,n) : \g \in \G_k, |x(\g,n)| > c\}}}
\newcommand{\xsc}{{\{s(\g,n) : \g \in \G_k, |x(\g,n)| > c\}}}
\newcommand{\pn}{{(\pi,n)}}
\newcommand{\pni}{{(\pi,n_i)}}
\newcommand{\pone}{{(\pi,1)}}
\newcommand{\ptwo}{{(\pi,2)}}

\newcommand{\cembeds}{\stackrel{c}{\hookrightarrow}}
\newcommand{\embeds}{\hookrightarrow}
\newcommand{\nem}{\not\embeds}
\newcommand{\da}{\downarrow}
\newcommand{\bs}{\backslash}
\newcommand{\we}{\wedge}
\newcommand{\nsim}{\not\sim}
\newcommand{\supp}{\operatorname{supp}}
\newcommand{\sgn}{\operatorname{sgn}}

\newcommand{\el}{\ell}
\newcommand{\G}{\Gamma}
\newcommand{\g}{\gamma}
\newcommand{\al}{\alpha}
\newcommand{\s}{\sigma}
\newcommand{\Si}{\Sigma}
\newcommand{\ep}{\epsilon}
\newcommand{\lam}{\lambda}
\newcommand{\Om}{\Omega}
\newcommand{\ap}{\aleph}
\newcommand{\D}{\Delta}
\newcommand{\sfm}{\mbox{\sf m}}

\newcommand{\onep}{{\frac{1}{p}}}

\newcommand{\hsp}{\hspace{1em}}

\newcommand{\nio}{^\infty_{n=0}}
\newcommand{\kio}{^\infty_{k=0}}
\newcommand{\nif}{^\infty_{n=1}}
\newcommand{\kif}{^\infty_{k=1}}
\newcommand{\ki}{^\infty_{k=n_0+n-1}}
\newcommand{\sil}{{\sum^l_{i=1}}}
\newcommand{\sF}{{\sum_{\pi\in F_i}}}
\newcommand{\sijl}{{\sum^l_{i=j+1}}}
\newcommand{\siji}{{\sum^\infty_{i=j+1}}}

\begin{document}

\title{Isomorphic classification of atomic weak $L^p$ spaces}
\author{Denny H. Leung}
\address{Department of Mathematics \\ National University of Singapore\\
         Singapore 0511}
\email{matlhh@@leonis.nus.sg}
\maketitle

\begin{abstract}
Let $\msp$ be a measure space and let $1 < p < \infty$.  The {\em weak
$L^p$}\/ space $\wlp$ consists of all measurable functions $f$ such
that 
\[
 \|f\| = \sup_{t>0}t^{\frac{1}{p}}f^*(t) < \infty,
\] 
where $f^*$ is the decreasing rearrangement of $|f|$.
It is a Banach space under a norm which is equivalent to the
expression above.  In this paper, we pursue the problem of 
classifying  weak $L^p$ spaces isomorphically when
$\msp$ is purely atomic.  It is also shown that if $\msp$ is a
countably generated $\sigma$-finite
 measure space, then $\wlp$ (if infinite dimensional) must
be isomorphic to either $\ell^\infty$ or $\seq$.
\end{abstract}



\section{Introduction}

In a series of papers \cite{L1,L2,L3}, the author studied the
isomorphic structures of the weak $L^p$ spaces defined on the measure
spaces $\N$, $[0,1]$, and $[0,\infty)$.  In this paper, we extend the
scope of this study to encompass all $\wlp$ spaces where $\msp$ is
purely atomic or countably generated and $\sigma$-finite.  
The main results obtained are:
\begin{enumerate}
\item The isomorphic classification of atomic weak $L^p$
spaces,
\item The fact that an infinite dimensional space $\wlp$, where $\msp$
is countably generated and $\sigma$-finite,
is isomorphic to either $\ell^\infty$ or $\seq$.
\end{enumerate}

Let us recall some standard definitions and set the notation.  Let
$\msp$ be an arbitrary measure space.  For $1 < p < \infty$, the {\em
weak $L^p$ space}\/ $\wlp$ is the space of all $\Sigma$-measurable
functions $f$ such that
\[ \|f\| = \sup_{t>0}t^{\frac{1}{p}}f^*(t) < \infty, \]
where $f^*$ denotes the decreasing rearrangement of $|f|$ \cite{LT2}.
One checks easily that 
\begin{equation}\label{norm}
 \|f\| = \sup_{c>0}c(\mu\{|f| > c\})^\onep , 
\end{equation}
and the latter expression is more convenient for computations.
It is well known
that $\|\cdot\|$ is equivalent to a norm under which $\wlp$ is a
Banach space.  If $\msp$ is a real interval $I$ endowed with Lebesgue
measure, we write $L^{p.\infty}(I)$; while $\seq$ and $\seq(m)$ will
stand for the weak $L^p$ spaces on $\N$ and $\{1,\ldots,m\}$
respectively, both with the counting measure.

The other notation follows mainly that of
\cite{LT1,LT2}. If $E$ and $F$ are Banach spaces, we write $E \sim F$
to indicate that they are isomorphic ($=$ linearly homeomorphic).
Then $E \embeds F$ (respectively, $E \cembeds F$) means that $E$ is 
isomorphic to a subspace (repectively, complemented subspace) of $F$.
For an arbitrary set $A$, we write $|A|$ for the cardinality of $A$,
and ${\cal P}(A)$ for the power set of $A$.  A measure space $\msp$ is
said to be {\em countably generated}\/ if there exists a sequence
$(A_n)\nif$ in $\Si$ such that the smallest $\s$-algebra generated by
the sets $(A_n)\nif$ and the $\mu$-null sets is $\Si$.  We will also
have occasion to use terms and notation concerning vector lattices,
for which the references are \cite{LT2,Sch}.  In
particular, two elements $a, b$ of a vector lattice are said to be {\em
disjoint}\/ if $|a| \we |b| = 0$.  A Banach lattice $E$
satisfies an {\em upper} $p$-{\em estimate}\/ if there is a constant
$M < \infty$ such that 
\[ \|\sum^n_{i=1}x_i\| \leq M(\sum^n_{i=1}\|x_i\|^p)^\onep \]
whenever $(x_i)^n_{i=1}$ is a pairwise disjoint sequence in $E$.
Using the expression (\ref{norm}), it is trivial to check that every
space $\wlp$ satisfies an upper $p$-estimate.

At this point, let us make 
some simple but useful observations. Let $\Gamma$ be an
arbitrary set, and let $w : \Gamma \to (0,\infty)$ be an arbitrary
function.  One can define a measure $\mu$ on the power set of $\Gamma$
by $\mu(A) = \sum_{\gamma\in A}w(\gamma)$ for all $A \subseteq
\Gamma$.  The resulting space $L^{p,\infty}(\G,{\cal P}(\G),\mu)$  
will be written as
$\seq(\G,w)$.  
If $w$ is identically $1$, we abbreviate it further to $\seq(\G)$.
It is clear that if $\msp$ is purely atomic, then
$\wlp$ is (linearly) isometric to $\seq(\G,w)$ for some pair $(\G,w)$.
The next observation further restricts the class of weight functions
$w$ which have to be considered.

\begin{pr}\label{scale}
Let $\G$ be an arbitrary set, and let $w_1, w_2 : \G \to (0,\infty)$
be such that $w_1(\g) \leq w_2(\g) \leq 2 w_1(\g)$ for all $\g \in
\G$.  Then $\seq(\G,w_1)$ and $\seq(\G,w_2)$ consist of the same
functions, and are isomorphic via the formal identity.
\end{pr}

It follows that it suffices to consider only weight functions $w$
which take values in the set $\set$. Such a $w$ will be called a {\em
modified weight function}\/ in the sequel. Another important
ingredient, used already in \cite{L3}, is Pelczynski's ``decomposition
method'' (see \cite[Proposition 21.5.12]{sem} and \cite[p.\ 54]{LT1}).
The {\em square}\/ of a Banach space $E$
is the Banach space $E \oplus E$.

\begin{thm}\label{pdm}
Let $E, F$ be Banach spaces such that
each is isomorphic to a complemented subspace of the
other.  Then they are isomorphic provided one of the following
conditions hold.
\begin{enumerate}
\item Both $E$ and $F$ are isomorphic to their squares,
\item $E$ is isomorphic to $(\sum\oplus E)_\elinf$.
\end{enumerate}
\end{thm}



\section{Countably many atoms}

In this section, we show that if $\G$ is countable, then $\seq(\G,w)$
is isomorphic to $\seq$ or $\ell^\infty$.  Together with the results
in \cite{L3}, this yields the isomorphic classification of $\wlp$ for
all countably generated measure spaces $\msp$.

\begin{pr}\label{Y}
Let $Y = (\sum\oplus\seq(n))_{\ell^\infty}$.  Then $\seq \sim Y$.
\end{pr}

\begin{pf}
Let $X = (\sum\oplus\seq(2^n))_{\ell^\infty}$.
It is known that 
\begin{equation}\label{X}
 \seq \sim \rline \cembeds X \cembeds \seq. 
\end{equation}
(See \S 2 of \cite{L3}, where $X$ is called $X_1$.)  Therefore,
\[ (\sum\oplus\seq)_{\ell^\infty} \cembeds (\sum\oplus X)_\elinf 
\cembeds
X \cembeds \seq . \]
Clearly $\seq \cembeds (\sum\oplus\seq)_{\ell^\infty}$ also.  Thus
$\seq \sim (\sum\oplus\seq)_{\elinf}$ by Theorem~\ref{pdm}.  It
follows immediately that $Y \cembeds \seq$.
Using (\ref{X})
again, we see that $\seq \cembeds X \cembeds Y$.  
Another application of Theorem~\ref{pdm} yields the desired result.
\end{pf}

The classification problem in the case of countably many atoms turns
out to hinge on two special cases.  We now describe the first case.

\begin{pr}\label{D}
Let $D = \{(i,j): i, j \in \N, i \geq j\}$, and let $(b_i)$ be a
strictly positive real sequence.  Define $\al : D \to (0,\infty)$ by
$\al(i,j) = b_i$.  If either $(i+1)b_{i+1} \leq b_i$ for all $i$, or
$b_{i+1} \geq (i+1)b_i$ for all $i$, then $\seq(D,\al) \sim \seq$.
\end{pr}

\begin{pf}
By Proposition~\ref{Y}, it suffices to show that $\seq(D,\al)$ is
isomorphic to $Y$.  Let $x \in \seq(D,\al)$.  For each $(i,j) \in D$,
let $y_i(j) = b_i^\onep x(i,j)$.  For each $i \in \N$, it is easy to
see that $y_i = (y_i(1),\ldots,y_i(i))$ is an element in $\seq(n)$ of
norm $\leq \|x\|$.  Hence the map $T: \seq(D,\al) \to Y$, $Tx = (y_i)$,
is a bounded linear map.  The proof will be complete if we can show
that $T$ is surjective.  Let $y = (y_i) \in Y$, and let $x : D \to \R$
be given by $x(i,j) = b_i^{-\onep}y_i(j)$.  If $\|x\| > 1$, there are
a $c > 0$ and a finite subset $B$ of $\{(i,j): |x(i,j)| > c\}$ 
such that $c\mu(B)^\onep > 1$, where $\mu$ is the measure on $D$
associated with the weight function $\al$.
Let $B_i = B \cap \{(i,j): 1 \leq j
\leq i\}$. Choose $i_1 \leq i_2$ such that $B_{i_1}, B_{i_2} \neq
\emptyset$, and $B_i = \emptyset$ if $i <
i_1$ or $i > i_2$.  Then, noting the fact that $|B_i| \leq i$,
\[ \mu(B) = \sum^{i_2}_{i=i_1}b_i|B_i| \leq \left\{
              \begin{array}{lll}
                3b_{i_1}|B_{i_1}| & \mbox{if}
                  & (i+1)b_{i+1} \leq b_i \hsp \mbox{for all $i$} \\
                3b_{i_2}|B_{i_2}| &  & 
                    b_{i+1} \geq b_i(i+1) \hsp \mbox{for all $i$}.
              \end{array}
             \right. \]
In either case, there is a $k \in \N$ such that $\mu(B) \leq
3b_k|B_k|$.  Thus $c(3b_k|B_k|)^\onep > 1$.  Therefore,
\[ cb_k^\onep|\{1 \leq j \leq k: |y_k(j)| > cb_k^\onep\}|^\onep \geq
cb_k^\onep|B_k|^\onep > 3^{-\onep} . \]
Hence $\|y\| > 3^{-\onep}$.  It follows that $x \in
\seq(D,\al)$ whenever $y \in Y$.  As $Tx = y$, $T$ is surjective, as
required. 
\end{pf}     
  
The second special case is the following.

\begin{pr}\label{Z}
Define $\s : \Z \to (0,\infty)$ by $\s(n) = 2^n$ for all $n \in \Z$.
Then $\seq(\Z,\s) \sim \elinf$.
\end{pr}

\begin{pf}
We claim that the map $T : \seq(\Z,\s) \to \elinf$, $T((a_n)) =
(2^{\frac{n}{p}}a_n)$, is an (onto) isomorphism. We continue to denote 
the
norm on $\seq(\Z,\s)$ by $\|\cdot\|$, while the norm on $\elinf$ will
be written as $\|\cdot\|_\infty$.  Also, we let $\mu$ be the measure
associated with the weight function $\s$.  Let $(a_n)$ be a normalized
element in $\seq(\Z,\s)$.  For any $j \in \Z$ such that $a_{j}
\neq 0$, and any positive $\ep < |a_{j}|$, 
\[ 2^{j} =
\mu\{j\} \leq \mu\{n : |a_n| > |a_{j}| - \ep\} \leq (|a_{j}|
- \ep)^{-p} . \]
Upon taking $\ep \da 0$, we see that $|a_{j}|2^{\frac{j}{p}} \leq
1$.  Therefore,  $T$ is bounded.  Conversely, assume that
$c = (c_n) \in \elinf$.  For any $\ep > 0$, 
\[ \mu\{n : \frac{|c_n|}{2^{\frac{n}{p}}} > \ep\} \leq \mu\{n :
\frac{\|c\|^p_\infty}{\ep^p} > 2^n\} = 2^{k+1} , \]
where $k$ is the largest integer such that $2^{k} <
\|c\|^p_\infty/\ep^p$.  Hence
\[ \ep(\mu\{n : \frac{|c_n|}{2^{\frac{n}{p}}} > \ep\})^\onep \leq \ep
2^{\frac{k+1}{p}} \leq 2^\onep\|c\|_\infty . \]
Therefore, $(2^{-\frac{n}{p}}c_n) \in \seq(\Z,\s)$.  Since 
$T((2^{-\frac{n}{p}}c_n)) = ((c_n))$, $T$ is onto. 
\end{pf}

We are now ready for the main result of this section.

\begin{thm}\label{count}
Let $\G$ be countably infinite, and let $w$ be a modified weight
function on $\G$.  Let $\G_n = \{\g \in \G : w(\g) = 2^n\}$ for all $n
\in \Z$.  Then
\begin{enumerate}
\item $\seq(\G,w) \sim \seq$ if and only if\/ $\sup_n|\G_n| = \infty$,
\item $\seq(\G,w) \sim \elinf$ if and only if\/ $\sup_n|\G_n| < \infty$.
\end{enumerate}
\end{thm}

\begin{pf}
Since $\G$ is countable, there is a collection $(A_\g)_{\g\in\G}$ of
pairwise disjoint Lebesgue measurable subsets of $[0,\infty)$ such
that $\lam(A_\g) = w(\g)$ for all $\g \in \G$, where $\lam$ is Lebesgue
measure. Let $\Si'$ be the $\s$-algebra generated by the sets $(A_\g)$.
Clearly, $\seq(\G,w) \sim L^{p,\infty}(\lam_{|\Si'})$, and the latter is
complemented in $\rline$ by a conditional expectation operator
\cite[Theorem 2.a.4]{LT2}.  Since $\rline \sim \seq$ \cite{L3},
$\seq(\G,w) \cembeds \seq$.  

Now assume that $\sup_n|\G_n| = \infty$.
We claim that $\seq \cembeds \seq(\G,w)$.
If $|\G_n| = \infty$ for some $n$, then obviously $\seq \sim
\seq(\G_n,w_{|\G_n}) \cembeds \seq(\G,w)$.
If $|\G_n|$ is finite for all $n$, choose a sequence
$(n_i)^\infty_{i=1}$ in $\Z$ 
such that $|\G_{n_i}| \geq i$ for all $i$, and either $2^{n_{i+1}} >
(i+1)2^{n_i}$ for all $i$, or $(i+1)2^{n_{i+1}} < 2^{n_i}$ for all $i$.
Taking $b_i$ in Proposition~\ref{D} to be $2^{n_i}$, we see that
\[ \seq \sim \seq(D,\al) \cembeds \seq(\bc_i\G_{n_i},w') \cembeds
\seq(\G,w) , \]
where $w'$ is the restriction of $w$ to $\bc_i\G_{n_i}$. Thus the
claim is proved.  In the proof of Proposition~\ref{Y}, it was observed
that $\seq \sim (\sum\oplus\seq)_\elinf$.  Hence we may apply
Theorem~\ref{pdm} to conclude that $\seq(\G,w)$ is isomorphic to
$\seq$. 

Next, consider the case when $\sup_n|\G_n| < \infty$.  Then
$\seq(\G,w) \cembeds \seq(\Z,\s) \oplus \cdots \oplus \seq(\Z,\s)$,
where the direct sum consists of a finite number of copies of
$\seq(\Z,\s)$. By Proposition~\ref{Z}, 
\[ \seq(\G,w) \cembeds \elinf \oplus \cdots \oplus \elinf \sim \elinf .
\]
Since it is easy to see  that $\elinf \cembeds \seq(\G,w)$,
$\seq(\G,w) \sim \elinf$ by Theorem~\ref{pdm}.  The converses follow
by simply observing that the spaces $\seq$ and $\elinf$ are not
isomorphic, since, for instance, the latter has the Dunford-Pettis
property while the former fails it.
\end{pf}

For the next theorem, we recall the following fact: if
$(\Om_i,\Si_i,\mu_i)$, $i = 1, 2$, are measure isomorphic measure
spaces, then $L^{p,\infty}(\Om_i,\Si_i,\mu_i)$, $i = 1, 2$, are
isometric \cite{Lotz}.  This follows easily from the observation that
the set of functions $f \in \wlp$ of the form 
\[ f = \bigvee_na_n\chi_{A_n}, \]
where $(a_n) \subseteq \R$, and $(A_n)$ is a pairwise disjoint
sequence in $\Si$, is dense in $\wlp$.

\begin{thm}\label{cogen}
Let $\msp$ be a countably generated $\sigma$-finite measure space.
Then $\wlp$ (if infinite dimensional) is isomorphic to $\seq$ or
$\elinf$.
\end{thm}
                                                    
\begin{pf}
Choose a sequence $(\Om_n)$ of pairwise disjoint $\Si$-measurable
sets 
of finite positive measure such that $\Om = \cup \Om_n$. 
Let $(\Om_n,\Si_n,\mu_n)$ be the restriction of the measure $\mu$ to
$\Om_n$. By a theorem of Carath\'{e}odory \cite[Theorem 15.2.2]{R}, 
each $(\Om_n,\Si_n,\mu_n)$ is measure isomorphic to a
measure subalgebra of the measure space $[0,\mu(\Om_n))$ with
Lesbesgue measure \cite{R}.  It follows readily that $\msp$ is measure
isomorphic to a measure subalgebra of $[0,\infty)$.
By \cite[Theorem 2.a.4]{LT2} and the main result in \cite{L3},
\[ \wlp \cembeds \rline \sim \seq. \]
For the converse, we may assume by the above argument that $\msp$ is a 
measure subalgebra of $[0,\infty)$.  In particular, this 
ensures the boundedness of conditional expectation operators.  Now
if $\msp$ is purely atomic, there is nothing to prove.  Otherwise,
$\msp$ has a measure subalgebra which is measure isomorphic to the
measure space $[0,1]$ with Lebesgue measure. Hence $\seq \sim \unit
\cembeds \wlp$. 
We may now apply Theorem~\ref{pdm} to conclude the proof.
\end{pf}

\noindent{\bf Remark}.  The $\sigma$-finiteness condition in Theorem 
\ref{cogen} cannot be removed, as the following example shows.  Let
$\msp$ be 
the space $\{-1,1\}^\N$ equipped with the usual product $\sigma$-algebra,
and the counting measure restricted to this $\sigma$-algebra.  
Then it is a countably generated measure space such that $\wlp$ is 
isometric to $\seq(\R)$.\\

We conclude this section with a result which verifies one of the
conditions for applying Theorem~\ref{pdm}.  This is useful in the next
section as it can be shown that $\seq(\G,w)$ is never isomorphic to
$(\sum\oplus\seq(\G,w))_\elinf$ when $\G$ is uncountable.

\begin{pr}\label{square}
Let $\G$ be an arbitrary infinite set, and let $w : \G \to
(0,\infty)$. Then $\seq(\G,w)$ is isomorphic to its square.
\end{pr}

\begin{pf}
It may be assumed that $w$ is a modified weight function.  For each $n
\in \Z$, let $\G_n = w^{-1}\{2^n\}$.  Let $N = \{n\in \Z: \G_n\,
\mbox{ is infinite}\}$, and $M = \Z \backslash N$.  
Let $\G' = \bc_{n\in N}\G_n$ and $\G'' = \bc_{n\in M}\G_n$.  
Clearly,  $\seq(\G',w_{|\G'})$ is isomorphic to its square.  
If $M$ is finite, then it
is easy to see that $\seq(\G,w) \sim \seq(\G',w_{|\G'})$, and the proof
is finished.
If $M$ is
infinite, $\seq(\G'',w_{|\G''}) \sim \seq$ or $\elinf$ by
Theorem~\ref{count}, and these spaces are isomorphic to their squares.
Hence $\seq(\G,w) \sim \seq(\G',w_{|\G'}) \oplus \seq(\G'',w_{|\G''})$
is also isomorphic to its square. 
\end{pf}


\section{The uncountable case. Preliminaries}

In this section, we present
some preliminary results concerning the isomorphic classification of 
$\seq(\G,w)$ spaces when $\G$ is uncountable.  Essentially,
isomorphically distinct $\seq(\G,w)$ spaces are distinguishable by
counting arguments.  The basic philosophy is that, first of all, the
larger the cardinality of the set $\G$, the ``bigger'' is the
resulting space.  Secondly, small atoms give rise to bigger spaces
than large atoms, since small atoms may be ``glued together'' to form
larger atoms. The precise formulation of the counting argument used is
given below.

\begin{pr}\label{disj}
Let $\G$ be an arbitrary set, and let $w$ be a weight function on
$\G$. Let $A \subseteq \seq(\G,w)$ and $\G_1 \subseteq \G$ be such
that $|A| > \max\{|\G_1|,\ap_0\}$.  Suppose also that there are
constants $K < \infty$, $r > 1$ such that 
\begin{equation}\label{dom}
\|\sum_{x\in F}\ep_xx\| \leq K|F|^{\frac{1}{r}} 
\end{equation}
for all finite subsets $F$ of $A$, and all $\ep_x = \pm 1$.  Then
there exists $A_1 \subseteq A$ such that $|A_1| = |A|$, and 
$x\chi_{\G_1}
= 0$ for all $x \in A_1$.
\end{pr}

\begin{pf}
Let $A_1 = \{x \in A : x\chi_{\G_1} = 0\}$.  For each $x \in A\bs A_1$,
there is a $\g \in \G_1$ such that $x(\g) \neq 0$.  Define a choice
function $f : A\bs A_1 \to \G_1$ such that $x(f(x)) \neq 0$ for all $x
\in A \bs A_1$.  If $|A_1| < |A|$, then $|A\bs A_1| = |A| > \ap_0$.  Hence
there exist $C \subseteq A\bs A_1$, and $n \in \N$ such that $|C| =
|A\bs A_1| = |A|$, and $|x(f(x))| \geq 1/n$ for all $x \in C$.  Now $|C|
= |A| > |\G_1| \geq |f(C)|$.  Therefore, there is a $\g_0 \in f(C)$ such
that $D = f^{-1}\{\g_0\} \cap C$ is infinite.  Note that $x \in D$ implies
$|x(\g_0)| \geq 1/n$.  Now for any finite subset $F$ of $D$, 
\[ \|\sum_{x\in F}\sgn x(\g_0)x\| \geq \sum_{x\in
F}|x(\g_0)|\|\chi_{\{\g_0\}}\| \geq \frac{|F|}{n}w(\g_0)^\onep . \]
As $D$ is infinite, this violates condition (\ref{dom}).
\end{pf}

\begin{pr}\label{size}
Let $w$ be a weight function on a set $\G$.  Let $A, B
\subseteq \seq(\G,w)$ be such that $|A| > \max\{|B|,\ap_0\}$.  Suppose
also that there are constants $K < \infty$, $r > 1$ such that
\[ \|\sum_{x\in F}\ep_xx\| \leq K|F|^{\frac{1}{r}} \]
for all finite subsets $F$ of $A$, and all $\ep_x = \pm 1$.  Then
there exists $C \subseteq A$, $|C| = |A|$, such that the elements of
$C$ are pairwise disjoint, and $|b| \we |c| = 0$ whenever $b \in B$, $c
\in C$. 
\end{pr}

\begin{pf}
For each $x \in \seq(\G,w)$, let $\supp x = \{\g \in \G : x(\g) \neq
0\}$. Clearly $|\supp x| \leq \ap_0$.  Therefore, $|\bc_{x\in B}\supp
x| \leq \max\{|B|,\ap_0\} < |A|$.   Let $\G_1 = \bc_{x\in B}\supp x$.
By Proposition~\ref{disj}, there is a subset $A_1$ of $A$, having the
same cardinality as $A$, such that $x\chi_{\G_1} = 0$ for all $x \in
A_1$. 
It remains to choose a pairwise
disjoint subset of $A_1$ of cardinality $|A|$.  This will be done by
induction. Choose $x_0$ arbitrarily in $A_1$.  Now suppose a pairwise
disjoint collection
$(x_\rho)_{\rho<\beta}$ has been chosen up to some ordinal $\beta <
|A| = |A_1|$. Since $|A|$ is a cardinal, $|\beta| < |A|$.  Hence
$|\bc_{\rho<\beta}\supp x_\rho| \leq \max\{|\beta|,\ap_0\} < |A|$.
Let $\G_2 = \bc_{\rho<\beta}\supp x_\rho$. Using
Proposition~\ref{disj} again, we find a $x_\beta \in A_1$ such that
$x_\beta\chi_{\G_2} = 0$.  It is clear that the collection
$(x_\rho)_{\rho\leq\beta}$ is pairwise disjoint.  This completes the
inductive argument.  Consequently, we obtain a pairwise disjoint
collection $C = (x_\rho)_{\rho<|A|}$ in $A_1$. As each $x \in C$ is
disjoint from each $b \in B$ by the first part, the proof is complete. 
\end{pf}
    
We first apply Proposition~\ref{size} to distinguish $\seq(\G,w)$
spaces according to the values of $p$ and $|\G|$.

\begin{thm}\label{diffG}
Let $(\G_i,w_i), i = 1,2$ be arbitrary, and let $1 < p, r < \infty$.  If
$|\G_1| \neq |\G_2|$, then $\seq(\G_1,w_1) \nsim \seqr(\G_2,w_2)$.
\end{thm}

\begin{pf}
If at least one of $\G_1, \G_2$ is finite, 
then the assertion is clear. Thus we may as well assume that $|\G_1| >
|\G_2| \geq \ap_0$.  Suppose $T : \seq(\G_1,w_1) \to \seqr(\G_2,w_2)$ is
an embedding.  For each $\g \in \G_1$, let $e_\g = \chi_{\{\g\}}$.
Define $x_\g = Te_\g/\|e_\g\|$ for all $\g \in \G_1$.  Let $A =
\{x_\g : \g \in \G_1\}$.  Then $|A| > \ap_0$.  Also, $\seq(\G_1,w_1)$
satisfies an upper $p$-estimate.  Therefore, there is a constant $K <
\infty$ such that 
\[ \|\sum_{\g\in F}\ep_\g\frac{e_\g}{\|e_\g\|}\| \leq K|F|^\onep \]
for all finite $F \subseteq \G_1$, and all $\ep_\g = \pm 1$.
Consequently, 
\[ \|\sum_{\g\in F}\ep_\g x_\g\| \leq K\|T\||F|^\onep . \]
Applying  Proposition~\ref{size} with $B = \emptyset$, we obtain a set
$C \subseteq A$, $|C| = |A|$, such that the elements of $C$ are
pairwise disjoint.  Thus $|C| > |\G_2|$, and $C$ is a  subset of 
$\seqr(\G_2,w_2)$ consisting of pairwise disjoint non-zero elements,
which is clearly impossible.
\end{pf}

\begin{thm}\label{diffp}
Let $\G$ be an arbitrary uncountable set, and let $1 < p \neq r <
\infty$. Then for any $w_i : \G \to (0,\infty)$, $i = 1,2$,
$\seq(\G,w_1) \not\sim \seqr(\G,w_2)$.
\end{thm} 

\begin{pf}
We may assume that $p > r$.  Suppose that $\seqr(\G,w_2) \embeds
\seq(\G,w_1)$. For any $n \in \Z$, let $\G_n = \{\g \in \G : 2^n \leq
w_2(\g) < 2^{n+1}\}$.  Then there is a $n_0$ such that $\G_{n_0}$ is
uncountable. Also, by Proposition~\ref{scale},
\[ \seqr(\G_{n_0}) \embeds \seqr(\G_{n_0},w_{2|\G_{n_0}}) \embeds
\seqr(\G,w_2) \embeds \seq(\G,w_1). \]
Now, let $T : \seqr(\G_{n_0}) \to \seq(\G,w_1)$ be an embedding, and
let $A = \{T\chi_{\{\g\}} : \g \in \G_{n_0}\}$.  By
Proposition~\ref{size}, there is a set $C \subseteq A$, $|C| = |A| =
|\G_{n_0}|$, such that the elements of $C$ are pairwise disjoint.
Since $T$ is an embedding, there exists $\ep > 0$ such that 
\[ \|\sum_{\g\in F}T\chi_{\{\g\}}\| \geq \ep|F|^{\frac{1}{r}} \]
for all finite sets $F$ such that
$\{\chi_{\{\g\}} : \g \in F\} \subseteq C$.  But this violates the
fact that 
$\seq(\G,w_1)$ satisfies an upper $p$-estimate.
\end{pf}

Before continuing, let us establish some more notation.  For the
remainder of the paper, we fix a $p$, $1 < p < \infty$. Let
$(\G_n)^\infty_{n=0}$ be an arbitrary collection of pairwise disjoint
sets, and define $w : \bc^\infty_{n=0}\G_n \to (0,\infty)$ by $w(\g)
= 2^n$ if $\g \in \G_n$. The resulting space $\seq(\bc\nio\G_n,w)$
will be 
denoted by $E(\G_n)$.  Similarly, we let $F(\G)$ be the space $\seq(\G
\times \N,w)$, where $\G$ is an arbitrary set, and $w(\g,n) =
2^{-n}$ for all $\g \in \G$, and $n \in \N$.

\begin{lem}\label{Gn}
Let $(\G_n)\nio$ be pairwise disjoint sets such that
$\max\{|\G_n|,\ap_0\} < \sup_m|\G_m|$ for all $n \geq 0$.  Then there
exists a sequence $(\Pi_n)\nio$ of pairwise disjoint sets such that
$|\Pi_n| \leq |\Pi_{n+1}|$,  
$\max\{|\Pi_n|,\ap_0\} < \sup_m|\Pi_m| = \sup_m|\G_m|$ 
for all $n \geq 0$, and $E(\G_n) \sim E(\Pi_n)$.
\end{lem}

\begin{pf}
>From the assumption, we see that for all $n \geq 0$, there is a $m >
n$ such that $|\G_m| > |\G_n|$. Let $n_0 = \min\{n \geq 0 : |\G_n|
\geq \ap_0\}$. Inductively, let $n_k = \min\{n > n_{k-1} : |\G_n| >
|\G_{n_{k-1}}|\}$ for $k \geq 1$. Now define
\[ \Pi_k = \left\{ \begin{array}{lcl}
                     \G_{n_0} \times \{k\} & \mbox{if} & 0 \leq k <
n_1 \\
                     \G_{n_j} \times \{k\} & & n_j \leq k < n_{j+1}, j
\geq 1 .
                   \end{array} 
           \right. \]
Then $(\Pi_k)\kio$ is a pairwise disjoint sequence of sets of
non-decreasing 
cardinalities such that
\[ \max\{|\Pi_n|,\ap_0\} < \sup_m|\Pi_m| = \sup_m|\G_m| \]  
for all $k$.
Since $|\G_n| \leq |\Pi_n|$ for all $n$,  it is clear
that $E(\G_n) \cembeds E(\Pi_n)$. On the other hand, as $\G_{n_j}$ is
infinite for all $j$, one can choose a sequence of pairwise disjoint
sets $(\D_n)\nio$ so that 
\[ |\D_n| = |\G_n| \hsp \mbox{ for all $n$,} \]
\[ \D_{n_0} = \G_{n_0} \times \{0,\ldots,n_1-1\} \times
\{1,\ldots,2^{n_1-n_0-1}\} , \]
and
\[ \D_{n_k} = \G_{n_k} \times \{n_k,\ldots,n_{k+1}-1\} \times
\{1,\ldots,2^{n_{k+1}-n_k-1}\}  \]
for $k \geq 0$.  The first condition assures that $E(\G_n)$ and
$E(\D_n)$ are isometric.  Define a map $s : \bc\nio\Pi_n \to
\bc\kio{\cal P}(\D_{n_k})$ as follows.  For $\pi \in \Pi_n, 0 \leq n
< n_0$, $\pi = (\g,n)$ for some $\g \in \G_{n_0}$.  In this case, let
$s(\pi) = \{(\g,n,1)\}$.  For $\pi \in \Pi_n$, $n_k \leq n <
n_{k+1}-1$, $k \geq 0$, $\pi = (\g,n)$ for some $\g \in \G_{n_k}$.
Then we let $s(\pi) = \{(\g,n,l) : 1 \leq l \leq 2^{n-n_k}\}$.  For $x
\in E(\Pi_n)$, let $Tx =
\bv\nio\bv_{\pi\in\Pi_n}x(\pi)\chi_{s(\pi)}$. 
We claim that $Tx \in E(\D_n)$.  Indeed, let $\mu$ and $\nu$ be the
measures associated with the spaces $E(\Pi_n)$ and $E(\D_n)$
respectively. Then for any $c > 0$,
\[ \mu\{|x| > c\} = \sum\nio 2^n|\{\pi\in\Pi_n : |x(\pi)| > c\}|. \]
While
\begin{eqnarray*}
\nu\{|Tx| > c\} & = & \sum\nio\sum_\xc\nu(s(\pi)) \\
  & = & \sum^{n_0-1}_{n=0}\sum_\xc\nu(s(\pi)) + 
        \sum^{\infty}_{n=n_0}\sum_\xc\nu(s(\pi)) \\
  & = & \sum^{n_0-1}_{n=0}\sum_\xc 2^{n_0} + 
        \sum\kio\sum^{n_{k+1}-1}_{n=n_k}\sum_\xc 2^{n-n_k}2^{n_k} \\
  & = & \sum^{n_0-1}_{n=0}2^{n_0}|\xc| + 
        \sum^{\infty}_{n=n_0}2^{n}|\xc| .
\end{eqnarray*}
It follows that
\[ \mu\{|x| > c\} \leq \nu\{|Tx| > c\} \leq 2^{n_0}\mu\{|x| > c\} . \]
Thus, not only does $Tx$ belong to $E(\D_n)$, $T$ is also an
isomorphic embedding.  Since $T(E(\Pi_n))$ is easily seen to be
complemented in $E(\D_n)$ by the conditional expectation with respect
to the $\s$-algebra generated by the sets $s(\pi)$, we get that
$E(\Pi_n) \cembeds E(\D_n) \sim E(\G_n)$.  Recalling
Proposition~\ref{square}, we apply Theorem~\ref{pdm} to finish the
proof.
\end{pf}         

We introduce one more class of spaces which is the analog of the class
$E(\G_n)$ for small atoms.  Let $(\G_n)\nif$ be a sequence of pairwise
disjoint sets, and let $w : \bc\nif\G_n \to (0,\infty)$ be the weight
function $w(\g) = 2^{-n}$ for $\g \in \G_n$.  The space
$\seq(\bc\nif\G_n,w)$ is denoted by $G(\G_n)$.
It will be shown presently, however, that this class
is isomorphically the same as the class $F(\G)$ introduced earlier.
The proof of the following lemma is essentially
the same as that of Lemma~\ref{Gn}. For a sequence of cardinals
$(\sfm_n)\nif$, we let $\limsup_n \sfm_n = \min_{n\geq 1}\sup_{k\geq
n}\sfm_k$.  

\begin{lem}\label{G}
Let $(\G_n)\nif$ be a sequence of pairwise disjoint sets such that
$\sup_m|\G_m| = \limsup_{m}|\G_m| > \ap_0$.  Then there
exists a sequence $(\Pi_n)\nif$ of pairwise disjoint sets such that
$|\Pi_n| \leq |\Pi_{n+1}|$ for all $n \geq 1$, 
$\sup_m|\Pi_m| = \sup_m|\G_m|$, and $G(\G_n) \sim G(\Pi_n)$.
\end{lem}

\begin{pr}\label{GF}
Let $(\G_n)\nif$ be a sequence of pairwise disjoint sets of
non-decreasing cardinalities such that $\sup_m|\G_m| > \ap_0$.
Then $G(\G_n) \sim F(\G)$, where $\G = \bc\nif\G_n$.
\end{pr}

\begin{pf}
Since $F(\G)$ is nothing but the space $G(\G'_n)$ where each $|\G'_n| =
|\G|$, we clearly have $G(\G_n) \cembeds F(\G)$. Choose $n_0$ such
that $\G_{n_0}$ is infinite.  For all $n \in \N$, let 
\[ \Pi_n =
(\bc^\infty_{k=n_0+n-1}\G_k) \times \{n\} .\]
Note that $|\Pi_n| = |\G|$ for all $n$.  Hence $F(\G) \sim G(\Pi_n)$.
Now let $\D_n = \G_n$ for $n < n_0$, and 
\[ \D_n = \G_n \times \{1,\ldots,n-n_0+1\} \times
\{1,\ldots,2^{n-n_0}\} \]
for $n \geq n_0$.  Then $|\D_n| = |\G_n|$ for all $n$.  Thus $G(\D_n)
\sim G(\G_n)$. We claim that $G(\Pi_n) \cembeds G(\D_n)$.  To this
end, define the subset $s(\g,n)$ of $\D_k$ by
\[ s(\g,n) = \{(\g,n,l) : 1 \leq l \leq 2^{k-n-n_0+1}\} \]
whenever $(\g,n) \in \G_k \times \{n\}$ for some $n \geq 1$, and $k
\geq n_0+n-1$. Then, for $x \in G(\Pi_n)$, let
\[ Tx = \bv\nif\bv\ki\bv_{\g\in\G_k}x(\g,n)\chi_{s(\g,n)} .\]
If $\mu$ and $\nu$ are the measures associated with the spaces
$G(\Pi_n)$ and $G(\D_n)$ respectively, and $c > 0$, then
\[ \mu\{|x| > c\} = \sum\nif\sum\ki 2^{-n}|\xgnc| . \]
While
\begin{eqnarray*}
 \nu\{|Tx| > c\} & = & \sum\nif\sum\ki\nu(\bc\xsc) \\
   & = & \sum\nif\sum\ki{\frac{2^{k-n-n_0+1}}{2^k}}|\xgnc| \\
   & = & 2^{1-n_0}\mu\{|x| > c\} .
\end{eqnarray*}
It follows that $T : G(\Pi_n) \to G(\D_n)$ is an isomorphic embedding.
Moreover, $T(G(\Pi_n))$ is complemented in $G(\D_n)$ by a conditional
expectation operator. This proves the claim.  The proposition follows
by applying Theorem~\ref{pdm}.
\end{pf}


\section{Isomorphic classes}

In this section, we show that if $\G$ is uncountable, then
$\seq(\G,w)$ falls into one of five isomorphic types. The fact that
these types are mutually exclusive will be dealt with in the next
section. To begin with, we consider the ``small'' and ``large'' atoms
separately. 

\begin{pr}\label{smallatoms}
Let $\G$ be an uncountable set, and let $w : \G \to \{2^{-n} : n \in
\N\}$ be a weight function.  Then $\seq(\G,w)$ is isomorphic to one of
the following:
\begin{enumerate}
\item $\seq(\G)$,
\item $\seq(\G) \oplus F(\Pi)$, where $\ap_0 < |\Pi| < |\G|$,
\item $F(\G)$.
\end{enumerate}
\end{pr}

\begin{pf}
Let $\G_n = w^{-1}\{2^{-n}\}$ for all $n \in \N$.
Choose $m_0$ such that $\sup_{n\geq m_0}|\G_n| = \limsup_{n}|\G_n|$, 
and denote this value by $\sfm$.  There is a strictly increasing
sequence $(m_k)\kio$, beginning with $m_0$, 
such that $(|\G_{m_k}|)\kio$ is non-decreasing, and
$\sup_{k\geq 0}|\G_{m_k}| = \sfm$. 
The
remainder of the proof is divided into cases.\\

\noindent\underline{Case 1}.  $\sfm \leq \ap_0$. \\
In this case, one of $\G_1,\ldots,\G_{m_0-1}$ has cardinality
$|\G|$.  It follows easily that
\[ G(\G_1,\ldots,\G_{m_0-1},\emptyset,\ldots) \sim \seq(\G_1) \oplus
\cdots \oplus \seq(\G_{m_0-1}) \sim \seq(\G) . \]
On the other hand, $\bc^\infty_{n=m_0}\G_n$ is finite or countably
infinite.  Hence
\[ G(\emptyset,\ldots,\emptyset,\G_{m_0},\G_{m_0+1},\ldots) \sim
\mbox{ $\seq$ or $\elinf$} \]
or is finite dimensional by Theorem~\ref{count}.  Since $G(\G_n)$ is
the direct  sum of the two
spaces above, we see that it is isomorphic to $\seq(\G)$. \\

\noindent\underline{Case 2}. $\ap_0 < \sfm < |\G|$. \\
Let the sequence $(m_k)\kio$ be as above.  Define 
\[ \Pi_n = \left\{ \begin{array}{lcll}
                    \emptyset & \mbox{ if } & n \neq m_k & \mbox{for
any $k$} \\ 
                    \G_{m_k} & & n = m_k & \mbox{for some $k$} .
                   \end{array}
           \right. \]
Then $G(\Pi_n) \cembeds \seq(\G,w)$.  By Lemma~\ref{G} and
Proposition~\ref{GF}, $G(\Pi_n) \sim F(\Pi)$, where $|\Pi| = \sfm$.
On the other hand, as $|\G_n| \leq |\Pi|$ for all $n \geq m_0$,
\[ F(\Pi) \sim G(\Pi_n) \cembeds
G(\emptyset,\ldots,\emptyset,\G_{m_0},\G_{m_0+1},\ldots) \cembeds 
G(\emptyset,\ldots,\emptyset,\Pi,\Pi,\ldots) \sim F(\Pi).\]
By Theorem~\ref{pdm}, we see that $F(\Pi) \sim
G(\emptyset,\ldots,\emptyset,\G_{m_0},\G_{m_0+1},\ldots)$.  Therefore,
\begin{eqnarray*} 
 \seq(\G,w) & \sim & G(\G_1,\ldots,\G_{m_0-1},\emptyset,\ldots) \oplus
  G(\emptyset,\ldots,\emptyset,\G_{m_0},\G_{m_0+1},\ldots) \\
  & \sim & \seq(\G_1\cup\cdots\cup\G_{m_0-1}) \oplus F(\Pi) \\
  & \sim & \seq(\G) \oplus F(\Pi) . 
\end{eqnarray*}
Note that the last isomorphism holds because
$|\G_1\cup\cdots\cup\G_{m_0-1}| = |\G|$ by the case assumption. \\

\noindent\underline{Case 3}. $\sfm = |\G|$. \\
Arguing as in Case 2, we see that 
\[ \seq(\G,w) \sim \seq(\G') \oplus F(\G) ,\]
where $\G' = \G_1\cup\cdots\cup\G_{m_0-1}$.  But it is clear that the
latter space is isomorphic to $F(\G)$.
\end{pf}

\begin{pr}\label{largeatoms}
Let $\G$ be an uncountable set, and let $w : \G \to \{2^n : n \in \N
\cup\{0\}\}$  be a weight function.  
Then $\seq(\G,w)$ is isomorphic to one of
the following: 
\begin{enumerate}
\item $\seq(\G)$,
\item $E(\Pi_n)$ for some sequence $(\Pi_n)\nio$ with non-decreasing
cardinalities such that $|\Pi_n| < \sup_m|\Pi_m| = |\G|$ for all $n
\geq 0$.
\end{enumerate}
\end{pr}

\begin{pf}
For all $n \geq 0$, let $\G_n = \{\g \in \G : w(\g) = 2^n\}$. Let $\D$
be the disjoint union of the sets $\G_n \times \{1,\ldots,2^n\}$, $n
\geq 0$. Then $|\D| = |\G|$.  For $\g \in \G_n$, let $s(\g) = \{(\g,l)
: 1 \leq l \leq 2^n\}$.  
Then let $Tx =
\bv\nio\bv_{\g\in\G_n}x(\g)\chi_{s(\g)}$ 
for $x \in \seq(\G,w)$. Computing as in the proof of Lemma~\ref{Gn},
one sees that $T$ is an isometric embedding of $\seq(\G,w)$ into
$\seq(\D) \sim \seq(\G)$.  Also, $T(\seq(\G,w))$ is complemented in
$\seq(\D)$ by a conditional expectation.  Hence $\seq(\G,w) \cembeds
\seq(\G)$. Now suppose there exists $n_0$ such that $|\G_{n_0}| =
\sup_m|\G_m| = |\G|$.  Then clearly $\seq(\G) \cembeds \seq(\G,w)$ as
well.  Therefore, $\seq(\G,w) \sim \seq(\G)$ by Theorem~\ref{pdm}.
Otherwise, $|\G_n| < \sup_m|\G_m|$ for all $n$. Since $\seq(\G,w) =
E(\G_n)$, we may apply Lemma~\ref{Gn} to complete the proof.
\end{pf}

\begin{thm}\label{class}
Let $\G$ be an uncountable set, and let $w : \G \to (0,\infty)$ be a
weight function.  Then $\seq(\G,w)$ is isomorphic to one of the
following types:
\begin{enumerate}
\item $\seq(\G)$,
\item \label{E} $E(\Pi_n)$ for some sequence $(\Pi_n)\nio$ with
non-decreasing 
cardinalities such that $|\Pi_n| < \sup_m|\Pi_m| = |\G|$ for all $n
\geq 0$,
\item $F(\G)$
\item $\seq(\G) \oplus F(\Pi)$, where $\ap_0 < |\Pi| < |\G|$,
\item $E(\Pi_n) \oplus F(\Pi)$, where $(\Pi_n)\nio$ satisfies the
condition in (\ref{E}), and $\ap_0 < |\Pi| < |\G|$.
\end{enumerate}
\end{thm}

\begin{pf}
We may assume that $w$ is a modified weight function.  For $n \in \Z$,
let $\G_n = w^{-1}\{2^n\}$.  Then let $\G_+ =
\bc\nio\G_n$, and $\G_- = \bc^{-1}_{n=-\infty}\G_n$, and let $u$
and $v$ be the restrictions of $w$ to these  respective sets.
If $|\G_+| \leq \ap_0$, then $\seq(\G_+,u) \cembeds \seq$ or $\elinf$
by Theorem~\ref{count}.  As it is easy to see, via Theorem~\ref{pdm},
that 
\[ \seq(\G_-,v) \sim \seq(\G_-,v) \oplus \seq \sim \seq(\G_-,v) \oplus
\elinf , \]
it follows that $\seq(\G,w) \sim \seq(\G_-,v)$ in this case.  By
Proposition~\ref{smallatoms}, $\seq(\G,w)$ is isomorphic to one of the
types (1), (3), or (4).  (Note that $|\G_-| = |\G|$ in this case.)
Similarly, if $|\G_-| \leq \ap_0$, then
$\seq(\G,w) \sim \seq(\G_+,u)$ is isomorphic to one of the types (1)
or (2) by Proposition~\ref{largeatoms}.  Finally, assume that both
$\G_-$ and $\G_+$ are uncountable. Combining
Propositions~\ref{smallatoms} and \ref{largeatoms}, we arrive at one
of the following possibilities.
\newcounter{icount}
\begin{list}%
{(\roman{icount})}
{\usecounter{icount}
 \setlength{\leftmargin}{9truemm}
 \setlength{\labelsep}{3truemm}}
\item $\seq(\G_-) \oplus \seq(\G_+)$,
\item $\seq(\G_-) \oplus E(\Pi_n)$, where $|\Pi_n| \leq |\Pi_{n+1}|$,
and $|\Pi_n| < \sup_m|\Pi_m| = |\G_+|$ for all $n$,
\item $\seq(\G_-) \oplus F(\Pi) \oplus \seq(\G_+)$, where $\ap_0 <
|\Pi| < |\G_-|$,
\item $\seq(\G_-) \oplus F(\Pi) \oplus E(\Pi_n)$, where $(\Pi_n)\nio$
and $\Pi$ satisfy the respective conditions in (ii) and (iii),
\item $F(\G_-) \oplus \seq(\G_+)$,
\item $F(\G_-) \oplus E(\Pi_n)$, where $(\Pi_n)\nio$ is as in (ii).
\end{list}
Clearly, (i) is isomorphic to $\seq(\G)$, which is case (1), and (iii)
is isomorphic to case (4). If $|\G_-| < |\G|$, then $|\G_+| = |\G|$.
Thus for case (ii), there exists $n_0$ such that $|\G_-| \leq
|\Pi_{n_0}|$.  Hence 
\[ \seq(\G_-) \oplus E(\Pi_n)
\cembeds E(\Pi_n) \cembeds \seq(\G_-) \oplus E(\Pi_n) . \]
Therefore, (ii) falls under case (2).  Similarly, (iv) belongs to the
class (5).  Also, (v) is the same as case (4), and (vi) is the same as
case (5) in this situation, since $\ap_0 < |\G_-| < |\G| = |\G_+|$.
It remains to consider the case when $\ap_0 < |\G_+| \leq
|\G| = |\G_-|$.
In this case, from the proof of Proposition~\ref{largeatoms}, we see
that 
\[ \seq(\G_-) \oplus E(\Pi_n) \cembeds \seq(\G_-) \oplus
\seq(\D) \sim \seq(\G), \]
where $\D$ is the disjoint union of the sets $\Pi_n \times
\{1,\ldots,2^n\}$. 
But 
\[ \seq(\G) \sim \seq(\G_-) \cembeds \seq(\G_-) \oplus E(\Pi_n) \]
as well.  Hence class (ii) falls under (1).  It follows that class
(iv) belongs to (4).  As for class (v),
\[ F(\G_-) \oplus \seq(\G_+) \sim F(\G) \oplus \seq(\G_+) \cembeds
F(\G) \]
since $|\G_+| \leq |\G|$. Since the reverse complemented inclusion is 
obvious, we see that (v) falls under (3).
Finally, using the identification of (ii) above,
\begin{eqnarray*}
 F(\G_-) \oplus E(\Pi_n) & \sim & F(\G_-) \oplus \seq(\G_-) \oplus
E(\Pi_n) \\
 & \sim & F(\G) \oplus \seq(\G) \\
 & \sim & F(\G) .
\end{eqnarray*}
Thus type (vi) is identified with (3).
\end{pf}


\section{Uniqueness}

The first half of this section will be devoted to proving that the
isomorphic classes in 
Theorem~\ref{class} are mutually exclusive.  Within each of the
classes (2), (4), and (5), there are different isomorphic types.  In
the latter half of the section, we identify  precisely the different
isomorphic types contained in each of these classes.

\subsection{}

The sole objective of this part is the

\begin{thm}\label{unique}
The isomorphic classes listed in Theorem~\ref{class} are mutually
exclusive. 
\end{thm}

The main steps are isolated in  the following lemmas. Throughout, $\G$
is a fixed uncountable set.

\begin{lem}\label{15}
$\seq(\G) \not\embeds E(\Pi_n) \oplus F(\Pi)$ for any sequence
$(\Pi_n)\nio$ of pairwise disjoint sets of non-decreasing
cardinalities such that $|\Pi_n| < \sup_m|\Pi_m| = |\G|$, and any set
$\Pi$ with $\ap_0 < |\Pi| < |\G|$.
\end{lem}

\begin{lem}\label{31}
$F(\Pi) \not\embeds \seq(\D)$ if $\Pi, \D$ are arbitrary sets such
that $\Pi$ is uncountable.
\end{lem}

\begin{lem}\label{34}
$F(\G) \not\embeds \seq(\G) \oplus F(\Pi)$ if $|\Pi| < |\G|$.
\end{lem}

Assuming these lemmas, we present the proof of Theorem~\ref{unique}.\\

\noindent{\em Proof of Theorem~\ref{unique}}.  
We adopt the notation $(a) \not\embeds (b)$, where $a, b \in
\{1,\ldots,5\}$, to indicate that no space of type (a) can be embedded
into a space of type (b).
By Lemma~\ref{15}, (1) $\nem$ (5).  Since any space of
type (2) is isomorphic to a subspace of a space of type (5), (1)
$\nem$ (2). Similarly, as a space of type (1) can be embedded into a
space of type (4), (4) $\nem$ (5).  From Lemma~\ref{31}, we see that
(3), (4), (5) $\nem$ (1).  Since a space of type (2) can be embedded
into a space of type (1) (see the proof of
Proposition~\ref{largeatoms}), it follows that (3), (4), (5) $\nem$
(2) as well.  By Lemma~\ref{34}, (3) $\nem$ (4).  Finally, as a space
of type (2) embeds into a space of type (1), a space of type (5)
embeds into a space of type (4).  Therefore, (3) $\nem$ (5).\hsp 
$\Box$\\

\noindent{\em Proof of Lemma~\ref{15}}.
For each $\g \in \G$, let $e_\g = \chi_{\{\g\}}$. 
Assume that $T : \seq(\G) \to E(\Pi_n) \oplus F(\Pi)$ is an embedding. 
Taking $A =
(Te_\g)_{\g\in\G}$ and $B = \emptyset$ in Proposition~\ref{size}, we
obtain a subset $\G'$ of $\G$ such that $|\G'| = |\G|$, and
$(Te_\g)_{\g\in\G'}$ is pairwise disjoint.  Now
\[ |\{\g \in \G' : \supp Te_\g \cap \Pi \neq \emptyset\}| \leq |\Pi| <
|\G| .\]
Discarding those elements that belong to the set above,we may assume
that 
$\supp Te_\g \subseteq \bc\nio\Pi_n$ for all $\g \in \G'$.  Since 
$T$ is an embedding, there is a $\ep > 0$ such that $\|Te_\g\| > \ep$
for all $\g \in \G'$.  Therefore, for every $\g \in \G'$, there is a
rational $c_\g > 0$ such that $c_\g(\nu\{|Te_\g| > c_\g\})^\onep > \ep$,
where $\nu$ is the measure associated with the space $E(\Pi_n)$.
Since $\G'$ is uncountable, there exists $\G'' \subseteq \G'$, $|\G''|
= |\G'| = |\G|$, and $c > 0$ such that $c(\nu\{|Te_\g| > c\})^\onep >
\ep$ for all $\g \in \G''$.  In particular, $\{|Te_\g| > c\} \neq
\emptyset$ for all $\g \in \G''$.  Let $n \in \N$ be given. Then
$|\Pi_1 \cup 
\cdots \cup \Pi_n| < |\G|$.  Using the disjointness of
$(Te_{\g})_{\g\in\G''}$, we obtain a $\g_0 \in \G''$ such that $\supp
Te_{\g_0} \subseteq \bc^\infty_{k=n+1}\Pi_k$.  
But then $\nu\{|Te_{\g_0}| > c\} \geq 2^{n+1}$.  Therefore,
\[ \|T\| \geq \|Te_{\g_0}\| \geq c(\nu\{|Te_{\g_0}| > c\})^\onep \geq
c2^{\frac{n+1}{p}} . \]
Since $n \in \N$ is arbitrary, we arrive at a contradiction.\hsp
$\Box$\\

\noindent{\em Proof of Lemma~\ref{31}}.  
Let $e_\pn$ be the characteristic function of $\{\pn\} \subseteq \Pi
\times \N$, and let $T : F(\Pi) \to \seq(\D)$ be an embedding.  Choose
$\ep > 0$ such that  $\|Te_\pn\| > \ep\|e_\pn\| = \ep 2^{-\frac{n}{p}}$
for all $\pn \in \Pi \times \N$.  Using the uncountability of $\Pi$,
we find a strictly positive 
sequence $(c_n)\nif$, and a sequence of
subsets $(\Pi_n)\nif$ of $\Pi$ such that, for each $n$,  $c^p_n$ is
rational, 
$|\Pi_n| = |\Pi|$, and 
\[ c_n(\nu\{|Te_\pn| > c_n\})^\onep > \ep 2^{-\frac{n}{p}} \]
for all $\pi \in \Pi_n$. Here $\nu$ is the measure associated with the
space $\seq(\D)$.  Apply Proposition~\ref{size} with $A =
\{Te_\pone : \pi \in \Pi_1\}$ and $B = \emptyset$ to obtain a {\em
countably}\/ infinite set $\Pi'_1 \subseteq \Pi_1$ such that
$\{Te_\pone : \pi \in \Pi'_1\}$ is pairwise disjoint.  Then apply the
proposition again with $A = \{Te_\ptwo : \pi \in \Pi_2\}$, $B =
\{Te_\pone : \pi \in \Pi'_1\}$ to obtain a countably infinite $\Pi'_2
\subseteq \Pi_2$ such that
\[ \{Te_\pone : \pi \in \Pi'_1\} \cup \{Te_\ptwo : \pi \in \Pi'_2\} \]
is pairwise disjoint.  Continuing inductively, we obtain, for each $n$,
a countably infinite $\Pi'_n \subseteq \Pi_n$ such that $\bc_n\{Te_\pn :
\pi \in \Pi'_n\}$ is pairwise disjoint. For all $\pi \in \Pi'_n$, $n
\in \N$, 
let $A_\pn = \{|Te_\pn| >
c_n\}$.  Then $\{A_\pn : \pi \in
\Pi'_n, n \in \N\}$ is a collection of pairwise disjoint 
subsets of $\D$ such that $\nu(A_\pn) > \ep^pc^{-p}_n2^{-n}$.  For the
remainder of the proof, we consider two cases.\\

\noindent\underline{Case 1}.  $\inf c_n > 0$. \\
Choose $\delta > 0$ such that $\inf c_n > \delta$.  For all $n$,
fix a  $\pi_n \in \Pi'_n$. As $\nu$ is just the counting measure on
$\D$, $\nu(A_{(\pi_n,n)}) \geq 1$ for all $n$.  Therefore, for any $k
\in \N$,
\[ \|\sum^k_{n=1}\frac{1}{c_n}Te_{(\pi_n,n)}\| \geq
\|\sum^k_{n=1}\chi_{A_{(\pi_n,n)}}\| \geq k^\onep . \]
Consequently,
\[ k^\onep \leq \|T\| \|\sum^k_{n=1}\frac{1}{c_n}e_{(\pi_n,n)}\| \leq
\frac{\|T\|}{\delta}\|\sum^k_{n=1}e_{(\pi_n,n)}\| \leq
\frac{\|T\|}{\delta} ,\]
a contradiction. \\

\noindent\underline{Case 2}. $\inf c_n = 0$. \\
Choose a strictly increasing sequence $(n_i)^\infty_{i=1}$ 
such that $c_{n_i} >
2c_{n_{i+1}}$ for all $i$.  Recall the fact that $c^p_n$ is rational
for all $n$.  Hence, for all $l \in \N$, there is a $m \in \N$ such
that the numbers $k_i = c^p_{n_i}2^{n_i}m$, $1 \leq i \leq l$, are all
integers. For $1 \leq i \leq l$, let $F_i$ be a subset of $\Pi'_{n_i}$
with $|F_i| = k_i$, and let $x = \sum^l_{i=1}c^{-1}_{n_i}\sF e_\pni$.
Then
\begin{eqnarray*}
 \|Tx\| & = & \|\sil\frac{1}{c_{n_i}}\sF Te_\pni\| \\
   & \geq & \|\sil\sF \chi_{A_\pni}\| \\
   & \geq & \{\sil\sF \nu(A_\pni)\}^\onep \\
   & \geq & \ep(ml)^\onep .
\end{eqnarray*}
Now let $c$ be a positive number.  If $c \geq 1/c_{n_l}$, then $\{|x|
> c\} = \emptyset$.  If $c^{-1}_{n_j} \leq c < c^{-1}_{n_{j+1}}$ for some
$1 \leq j < l$, then
\[ \mu\{|x| > c\} = \sijl\frac{k_i}{2^{n_i}} , \]
where $\mu$ is the measure associated with $F(\Pi)$.
Therefore,
\begin{eqnarray*}
 c(\mu\{|x| > c\})^\onep & \leq &
  \frac{1}{c_{n_{j+1}}}(\sijl\frac{k_i}{2^{n_i}})^\onep \\
  & = & \frac{m^\onep}{c_{n_{j+1}}}(\sijl c^p_{n_i})^\onep \\
  & = & m^\onep(\siji(\frac{c_{n_i}}{c_{n_{j+1}}})^p)^\onep \\
  & \leq & m^\onep(\siji(\frac{1}{2^{i-j-1}})^p)^\onep \\
  & \leq & 2m^\onep .
\end{eqnarray*}
The same computation holds when $c < 1/c_{n_1}$ if we replace $j$ by
$0$.  Hence $\|x\| \leq 2m^\onep$.
Thus
\[ \ep(ml)^\onep \leq \|Tx\| \leq \|T\|\|x\| \leq \|T\|2m^\onep . \]
We have arrived at yet another contradiction since $l \in \N$ is
arbitrary.\hsp
$\Box$\\

\noindent{\em Proof of Lemma~\ref{34}}.
As before, let $e_{(\g,n)}$ be the characteristic function of
$\{(\g,n)\} \subseteq \G \times \N$.  Similiarly define $e_\pn$ for
$\pn \in \Pi \times \N$. Suppose there is an embedding $T
: F(\G) \to \seq(\G) \oplus F(\Pi)$. For each $n \in \N$, let 
\[ A_n =
\{Te_{(\g,n)} : (\g,n) \in \G \times \N\}. \]
Also, let $B = \{e_\pn : \pn \in \Pi \times \N\}$.  Apply
Proposition~\ref{size} to obtain sets $\G_n \subseteq \G$, such that
$|\G_n| = |\G|$, and $\{Te_{(\g,n)} : \g \in \G_n\} \cup B$ is
pairwise disjoint for each $n$.  Then $Te_{(\g,n)} \in \seq(\G)$ for
all $\g \in \G_n$.  
A contradiction may now be derived as in the proof of Lemma~\ref{31}.\hsp
$\Box$

\subsection{}

In this last part, we derive the conditions which allow us to
distinguish between isomorphically distinct spaces that lie in the
same class (2), (4) or (5).  It was observed very early on
(Proposition~\ref{scale}) that scaling all atoms by factors bounded
away from $0$ and $\infty$ does not alter the isomorphic class of the
space.  We will see that within the same type (2), (4), or (5), this
is all one can do to stay in the same isomorphic class.

Once again, fix an uncountable set $\G$. A sequence of sets
$(\Pi_n)\nio$ will be called {\em admissible}\/ if it satisfies the
condition in
(2) of Theorem~\ref{class}.  We begin by considering type
(5) spaces.

\begin{thm}\label{type5}
Let $(\Pi_n)\nio$, $(\D_n)\nio$ be admissible sequences of sets, and
let $\Pi$, $\D$, be sets such that $\ap_0 < |\Pi|, |\D| < |\G|$.  Then
$E(\Pi_n) \oplus F(\Pi) \sim E(\D_n) \oplus F(\D)$ if and only if
$|\Pi| = |\D|$, and there exists $k \in \N$ such that 
\begin{equation}\label{PD}
 |\Pi_n| \leq |\D_{n+k}| \leq |\Pi_{n+2k}| 
\end{equation}
for all $n \in \N \cup \{0\}$.
\end{thm}

\begin{pf}
Clearly $F(\Pi) \sim F(\D)$ if $|\Pi| = |\D|$.  Now assume the
existence of $k \in \N$ such that $|\Pi_n| \leq |\D_{n+k}|$ for all $n
\in \N \cup \{0\}$. Choose a map $j : \bc\nio\Pi_n \to \bc\nio\D_n$ so
that, for each $n$, $j_{|\Pi_n}$ maps $\Pi_n$ injectively into
$\D_{n+k}$.  Define $T : E(\Pi_n) \to E(\D_n)$ by
\[ Tx = \bigvee\{x(\pi)e_{j(\pi)} : \pi \in \bc\nio\Pi_n\} . \]
Using calculations similar to those which appeared before, one sees
that $T$ maps $E(\Pi_n)$ isomorphically onto a complemented subspace
of $E(\D_n)$. 
Therefore, condition (\ref{PD}) yields the two way embedding
\[ E(\Pi_n) \cembeds E(\D_n) \cembeds E(\Pi_n) . \]
Consequently, $E(\Pi_n) \sim E(\D_n)$ by Theorem~\ref{pdm}.  This
proves the ``if'' part.

Conversely, assume that $E(\Pi_n) \oplus F(\Pi) \sim E(\D_n) \oplus
F(\D)$.  In particular, $F(\D) \embeds E(\Pi_n) \oplus F(\Pi) \embeds
\seq(\G) \oplus F(\Pi)$.  It follows from Lemma~\ref{34} that $|\D|
\leq |\Pi|$.  By symmetry, $|\Pi| = |\D|$.  Now 
\[ |\Pi| = |\D| < |\G| = \sup_m|\D_m| = \sup_m|\Pi_m| .\]
Thus there exists $n_0$ such that $|\Pi_{n_0}|, |\D_{n_0}| > |\Pi| =
|\D|$.  Since $E(\Pi_n) \sim E(\Pi_{n_0+n})$, and $E(\D_n) \sim
E(\D_{n_0+n})$, we may assume that $|\Pi_0|, |\D_0| > |\Pi| = |\D|$.
Now suppose that for no $k \in \N$ does $|\Pi_n| \leq |\D_{n+k}|$ hold
for all $n$.  Then there is a strictly increasing sequence $(n_k)\kif$
in $\N$ such that $|\Pi_{n_k}| > \max\{|\Pi_{n_{k-1}}|,
|\D_{n_k+k}|\}$ for all $k \in \N$.  Let $T : E(\Pi_n) \to E(\D_n)
\oplus F(\D)$ be an embedding, and let $e_\pi = \chi_{\{\pi\}}$ for all
$\pi \in \bc\nio\Pi_n$.  
Similarly define $e_{(\delta,n)}$ for $(\delta,n) \in \D \times \N$,
and $e_\delta$ for $\delta \in \bc\nif\D_n$. 
For each $k \in \N$, let $A = \{Te_\g : \g \in
\Pi_{n_{k+1}}\}$, and 
\[ B = \{e_{(\delta,n)} : (\delta,n) \in \D \times \N\} \cup
       \{Te_\pi : \pi \in \bc^k_{i=1}\Pi_{n_i}\} \cup
       \{e_\delta : \delta \in \bc^{n_{k+1}+k+1}_{j=1}\D_j\} . \]
Then
\[ |B| = |\D| + \sum^k_{i=1}|\Pi_{n_i}| +
\sum^{n_{k+1}+k+1}_{j=1}|\D_j| < |\Pi_{n_{k+1}}| = |A| . \]
Using Proposition~\ref{size}, one can choose $\Pi'_{k+1} \subseteq
\Pi_{n_{k+1}}$ such that $|\Pi'_{k+1}| = |\Pi_{n_{k+1}}|$, and $B \cup
\{Te_\pi : \pi \in \Pi'_{k+1}\}$ is pairwise disjoint. 
Then $\{Te_\pi : \pi \in \bc\kif\Pi'_{k+1}\}$ is a pairwise disjoint
subset of $E(\D_n)$ such that $\supp Te_\pi \subseteq
\bc^\infty_{j=n_{k+1}+k+2}\D_j$ whenever $\pi \in \Pi'_{k+1}$.
Let $\nu$ be the measure associated with the space $E(\D_n)$.
Arguing as in the proof of Lemma~\ref{31}, one obtains $\ep > 0$, and,
for each $k \in \N$, $c_{k+1} > 0$, and an infinite subset
$\Pi''_{k+1}$ of $\Pi'_{k+1}$, such that
\[ c_{k+1}(\nu\{|Te_\pi| > c_{k+1}\})^\onep > \ep\|e_\pi\| = \ep
2^{\frac{n_{k+1}}{p}} \]
for all $\pi \in \Pi''_{k+1}$.  If $\inf c_{k+1} \neq 0$, choose $\eta
> 0$ such that $\inf c_{k+1} > \eta$. Then, for all $k \in \N$, and
all $\pi \in \Pi''_{k+1}$, 
\[ \{|Te_\pi| > \eta\} \supseteq \{|Te_\pi| > c_{k+1}\} \neq \emptyset .
\]
Since the first set is contained in 
$\bc^\infty_{j=n_{k+1}+k+2}\D_j$, we
see that $\nu\{|Te_\pi| > \eta\} \geq 2^{n_{k+1}+k+2}$.  Therefore,
\[ 2^{\frac{n_{k+1}}{p}}\|T\| = \|T\|\|e_\pi\| \geq \|Te_\pi\| \geq \eta
2^{\frac{n_{k+1}+k+2}{p}}. \]
This is a contradiction since $k$ is arbitrary.  If $\inf c_{k+1} =
0$, we argue as in the corresponding Case 2 in  the proof of
Lemma~\ref{31} to obtain a contradiction.  The only difference in the
present situation is that $m$ is to be chosen so that the numbers
$2^{-n_{k_i}}c^p_{k_i}m$ are integers.
With this, we have shown that $E(\Pi_n) \not\embeds E(\D_n) \oplus
F(\D)$ if there is no $k$ such that $|\Pi_n| \leq |\D_{n+k}|$ for all
$n$. It follows easily that  condition (\ref{PD}) holds whenever
$E(\Pi_n) \oplus F(\Pi)$ and $E(\D_n) \oplus F(\D)$
are isomorphic.
\end{pf}

Contained in the above proof is the 

\begin{thm}
Let $(\Pi_n)\nio$ and $(\D_n)\nio$ be admissible sequences of sets.
Then $E(\Pi_n)$ and $E(\D_n)$ are isomorphic
if and only if condition (\ref{PD}) of
Theorem~\ref{type5} holds.
\end{thm}

It remains to examine type (4) spaces.

\begin{thm}
Let $\Pi, \D$ be sets such that $\ap_0 < |\Pi|, |\D| < |\G|$.  Then
$\seq(\G) \oplus F(\Pi) \sim \seq(\G) \oplus F(\D)$ if and only if
$|\Pi| = |\D|$.
\end{thm}

\begin{pf}
The ``if'' part is trivial. Suppose $|\Pi| > |\D|$, and $T : F(\Pi)
\to \seq(\G) \oplus F(\D)$ is an embedding. Let $e_{(\pi,n)} =
\chi_{\{(\pi,n)\}}$ for all $(\pi,n) \in \Pi \times \N$.  Similarly
for $e_{(\delta,n)}$, $(\delta,n) \in \D \times \N$.
Let $A = \{Te_{(\pi,1)} : \pi \in \Pi\}$ and let $B =
\{e_{(\delta,n)} : (\delta,n) \in \D \times \N\}$.  
Apply Proposition~\ref{size} to obtain a countably infinite
$\Pi_1 \subseteq \Pi$ such that $\{Te_{(\pi,1)} : \pi \in \Pi_1\} \cup
B$ is pairwise disjoint.  Continue inductively as in the proof of
Lemma~\ref{31} to obtain a sequence
$(\Pi_n)\nif$ of countably infinite subsets of $\Pi$ such that
\[ \{Te_{(\pi,n)} : \pi \in \Pi_n, n \in \N\} \cup B \]
is pairwise disjoint.  In particular, $Te_{(\pi,n)} \in \seq(\G)$ for
all $\pi \in \Pi_n$, $n \in \N$. We may now proceed as in the proof of
Lemma~\ref{31} to obtain a contradiction.
\end{pf}

\noindent{\bf Remark}.  The only remaining question in the isomorphic 
classification of atomic weak $L^p$ spaces is whether the spaces $\seq$
and $\seqr$ are isomorphic if $1 < p \neq r < \infty$.  The techniques 
of the present paper, which distinguish Banach spaces by showing that
one does not embed into the other, cannot be applied to this question.
This is because of the simple observation that every space $\seq$ 
contains an isomorphic copy of $\elinf$, and is also isomorphic to a 
subspace of $\elinf$.  The author has been able to obtain the partial 
result that $\seq$ and $\seqr$ are not isomorphic if $1 < p < r < \infty$
and $p \leq 2$.



\bibliographystyle{standard}
\bibliography{tref}

\end{document}